\documentclass[12pt]{article}
\input{epsf}

\usepackage{amsmath,amssymb,amsfonts,graphicx}
\newtheorem{thm}{Theorem}
\newtheorem{lemma}{Lemma}

\newtheorem{cor}{Corollary}

\newtheorem{rem}{Remark}

\begin{document}

\begin{center}
{\Large \bf  A Proof of the Erd\'os-Ulam Problem Assuming}\\
\medskip
 {\Large \bf Bomberi-Lang Conjecture}\\
\vspace*{0.7cm}
{\bf  Jafar Shaffaf}\\

\end{center}

\vspace*{1cm}

\begin{abstract}
\noindent A rational set in the plane is a point set with all its pairwise distances rational. Ulam and Erd\'os conjectured in 1945 that there is no dense rational set in the plane.  In this paper we associate special surfaces, we call them distance surfaces, to finite (rational) sets in the plane.   We prove that under a mild condition on the points in the rational set $S$ the associated distance surface in $\mathbb{P}^3$ is a surface of general type.  Also if we assume Bombieri-Lang  Conjecture in arithmetic algebraic geometry we can answer the Erd\'os-Ulam problem.
\end{abstract}

\begin{center}
\end{center}

\section{Introduction}

A rational (respectively integral) distance set is a subset $S$ of the plane $\mathbb{R}^2$ such that for all $s,t \in S$, the distance between $s$ and $t$ is a rational number (respectively, is an integer). On any line, one can easily find an infinite rational set that is in fact dense.  One can also easily find an everywhere dense rational subset of the unit circle. However it is not known if there is a rational set with $8$ points in general position, that is, no $3$ of them on a line and no $4$ on a circle.
In 1945, Anning and Erd\'os \cite{AnningErdos}  proved that any infinite integral set must be contained in a line. In 1945, Ulam conjectured that there is no everywhere dense rational set in the plane, see Problem III.5 in \cite{Guy}.  A deep conjecture of Erd\'os \cite{Erdos} states that an infinite rational set has to be very restricted.

\medskip

Huff \cite{Huff} considered rational distance sets $S$ of the following form: given distinct $a,b \in  \mathbb{Q}^{*}$, $S$ contains the four points $(0, \pm a)$ and $(0, \pm b)$ on the $y$-axis, plus points $(x, 0)$ on the $x$-axis, for some $x \in \mathbb{Q}^{*}$.  Such a point $(x, 0)$ must then satisfy the equations $x^2 + a^2 = u^2$ and $x^2+b^2 = v^2$ with $u, v \in \mathbb{Q}$. The system of associated homogeneous equations $x^2+a^2z^2 = u^2$ and $x^2 + b^2z^2 = v^2$ defines a curve $C(a^2, b^2)$ of genus $1$ in $\mathbb{P}^3$. Huff \cite{Huff} and  Peeples \cite{Peeples} provided examples in which the elliptic curve $ C(a^2,b^2)$ has a positive rank over $\mathbb{Q}$, thus exhibiting examples of infinite rational distance sets that are not contained in a line or in a circle. These remain to this day the 'largest' known such examples.  Assuming Bombieri-Lang conjecture, in this paper we show that this is actually the largest possible example (see Corollary \ref{cor1}). 

\medskip

Using Mordell conjecture, proved by Faltings, J. Solymosi and F. De Zeeuw \cite{SZ} proved  that lines and circles are the only irreducible algebraic curves that contain an infinite rational set. They also showed that if a rational set S has infinitely many points on a line or on a circle, then all but 4 (respectively 3) points of S are on the line (respectively circle). This answers questions of Guy, Problem D20 in \cite{Guy}, and Pach, Section 5.11 in \cite{BassMoserPach}.  For a rational distance subset of the plane intersecting any line in finitely many points,  it was shown in \cite{MS} that there is a uniform bound on the number of these intersections.  The main tool in \cite{MS} is a conditional uniform boundedness theorem on the number of rational points of algebraic curves of genus $ g \ge  2$ on a fixed number field $K$ (\cite{CHM1}). 

\medskip

For one of our main results we need to assume the Bombieri-Lang Conjecture.
The Lang Conjecture generalize the Mordell Conjecture to varieties of higher dimension.  For a more complete discussion of Lang Conjectures see for example \cite{Lang}.
Analogous of the Mordell Conjecture for varieties of {\it general type} - a natural generalization of the notion of curves of genus $g \ge 2$ to higher dimensions - is as follows (for the  precise definition see for instance \cite{Lang}):


\medskip

\noindent 
{\bf Weak Lang Conjecture. }  Let $X$ be a projective variety of general type defined over a number field $K$.  Then the set of rational points $X(K)$ is not Zariski dense in $X$.

It is worth mentioning one of the main consequences of the Weak Lang Conjecture in number theory and arithmetic geometry.

\medskip

\noindent
{\bf Uniform Bound Theorem.}	The weak Lang conjecture implies that, for every number field $K$ and for every $g \ge  2$, there exists a number $B(K,g)$ such that no curve of genus $g$ defined over $K$ has more than $B(K, g)$ points defined over $K$.

\medskip

For the main references of the Lang Conjectures and its consequences about the distribution of rational points on curves see \cite{CHM1}, \cite{CHM2}, \cite{Lang} and \cite{Ca2}.

\medskip

The result of this paper should be considered as  a witness for the power of  the Bombieri-Lang conjecture, rather than a real proof for this problem.  However, to the best of our knowledge we are far beyond a proof of the Lang conjecture compared with the same conjecture for curves that is,  Faltings Theorem (Mordell conjecture).

\medskip

In this paper we associate a surface to a finite (rational) subset \\ $A =
\{(\alpha_1, \beta_1), \ldots, (\alpha_m, \beta_m)\} \subset \mathbb{R}^2$ as follows:

\begin{equation}
\label{eq1}
z^2 = \prod_{i=1}^{m}\{(x-\alpha_i)^2 + (y-\beta_i)^2\},
\end{equation}
and we call it distance surface.  Also we will see that the points in the rational distance set play the role of the rational points on these distance surfaces over an appropriate number field.   \medskip

The following theorem is one of our main result (see Section 3):

\begin{thm}
\label{gentype}
Let $m=2g+2$ be an even integer number with $g \ge 2$ and assume $A=\{(\alpha_1 , \beta_1) , \ldots , (\alpha_m, \beta_m)\} $ be a finite subset of points in the plane not all on a line. Then the associated distance surface $S$ (the projectivization of the equation {\rm (1.1) } in $\mathbb{P}^3$) is of general type.
\end{thm}
The above theorem is our key tool to deduce the following main result of this paper.
\begin{thm}
\label{ulam}
Assuming Bombieri-Lang Conjecture, there is no dense rational distance subset in the plane.
\end{thm}
In the end, as a corollary of the proof of the above theorem we show that an infinite rational set in the plane should be very restricted.
\begin{cor}
\label{cor1}
Let $S$ be a rational set with infinitely many points not all on a line. Then all but at most $4$ ($3$) points lie on a line (circle).
\end{cor}

\medskip



The paper is organized as follows.
 In Section 2 we recall some preliminaries and known results about rational distance sets. 
 In Section 3 we explain the construction of a special surface in $\mathbb{P}^3$, using a rational distance set for which the points in the rational set play the role of its rational points over an appropriate number field and we also prove Theorem 1.1. At the end, assuming the Bombieri-Lang conjecture in arithmetic algebraic geometry we answer the Erd\'os-Ulam problem.

\section{Preliminaries and Known Results on Rational Distance Sets}
Rationality of distances in $\mathbb{R}^2$ is clearly preserved by translations, rotations, and uniform scaling ($(x, y) \rightarrow (\lambda x , \lambda y)$ with $\lambda \in \mathbb{Q}$).  We call a transformation $T: \mathbb{R}^2 \rightarrow \mathbb{R}^2$ a similarity transformation if it can be written as the composition of  translations, rotations and uniform scaling. Also rational sets are preserved under certain central inversions (an inversion with respect to a point in the rational set and with a rational radius), more precisely we quote the following lemma from \cite{SZ}.\\

\begin{lemma}
If we apply inversion to a rational set $S$, with the center $x \in S$ and a rational radius, then the image of $S \setminus \{x\}$ is a rational set.
\end{lemma} 

A priori, points in a rational set have arbitrary coordinates. However, after moving two of the points to two fixed rational points by translating, rotating and scaling, the points are almost rational points. The following simple lemma is well known (see \cite{K}).
\begin{lemma}
\label{a0}
 For any rational set $S$ there is a square free integer $k$ such that if a similarity transformation $T$ transforms two points of $S$ into $(0,0)$ and $(1,0)$, then any point in $T(S)$ is of the form $(r_1, r_2\sqrt{k}) $ , $r_1, r_2 \in \mathbb{Q}$. 
\end{lemma}

{\bf Proof }
Let $(0,0)$ and $(1,0)$ and $(x, y)$ be members of the rational set  $S$.  Then we have $x^2 + y^2 = q_1^2$ and $(x-1)^2 + y^2 = q_2^2$ such that $q_1 , q_2 \in \mathbb{Q}$.   By solving these two equations we have $x= \frac{q_1^2 - q_2^2 +1}{2} $, and so $x$ is a rational number.   By substitution of the value of $x$ and a simple manipulation for $y$ it can be easily shown that  $y= r\sqrt{k}$ where $r \in \mathbb{Q}$ and $k$ is square free.  For uniqueness  of $k$ let $p_1=(r_1, r_2\sqrt{k})$ and $p_2=(r_1', r_2'\sqrt{k'})$ be in the rational set $S$.  Hence the distance between $p_1$ and $p_2$ is a rational number.  By writing the square of the distance between two points $p_1$and $p_2$ it follows that the number $2r_2r_2' \sqrt{kk'}$ should be a rational number and hence $\sqrt{kk'}$ is rational.  Since $k$ and $k'$ are square free therefore we have $k=k'$.  $\blacksquare$\\
 
{

Note that we have the following two important theorems from \cite{SZ} :\\

\begin{thm}({\rm \cite{SZ}})Every rational subset of the plane has only  finitely many points in common with an algebraic curve defined over $\mathbb{R}$, unless the curve has a  line or circle as a component.
\end{thm}

\begin{thm}({\rm \cite{SZ}})
If a rational set $S$ has infinitely many points on a line (respectively circle) then all but $4$ (respectively $3$) points of $S$ are on the line (respectively circle).
\end{thm}
As we see the previous two theorems give a complete description of rational distance sets with infinitely many points on a line.  Hence if the Erd\'os-Ulam conjecture is false; i.e. there is a dense rational set in $\mathbb{R}^2$, then Theorem 4 implies that every line intersects $S$ in finitely many points.   
The following main result in \cite{MS} gives a uniform bound on the cardinality of these intersections.
\begin{thm}
\label{a1}
Assume Uniform Bound Conjecture and let $S$ be an infinite rational subset of the plane which has finite intersection with each line.  Then there exist an integer number $B(S)$ such that for each line $L$ in the plane the intersection $S \cap L$ has at most $B(S)$ points.  Moreover $B(S)\leq B(\mathbb{Q}, 2)$, in particular there is an upper bound for $B(S)$ independent of $S$.
\end{thm}

\section{Proofs}
Let $A=\{(\alpha_1 , \beta_1) , \cdots , (\alpha_m, \beta_m)\}$ be a subset of $\mathbb{R}^2$ not all on a line  and assume $m$ is an even integer number of the form $m=2g+2$ with $g \ge 2$.  Recall that in the introduction we associate to the subset $A$ the following surface $X$ 

\begin{equation}
\label{distanc}
z^2 = \prod_{i=1}^{m}\{(x-\alpha_i)^2 + (y-\beta_i)^2\},
\end{equation}
and we called it distance surface.  For finding the singularities of this surface we can rewrite the equation of the surface as follows :

\begin{align*}
z^2= \prod_{j=1}^{m}\{(x-\alpha_j)+i(y-\beta_j)\}\{(x-\alpha_j)-i(y-\beta_j)\},
\end{align*}
or equivalently we can write

\begin{align*}
z^2= \prod_{j=1}^{m} (x+iy-(\alpha_j + i \beta_j))(x-i y -(\alpha_j -i \beta_j)).
\end{align*}
Now by making the change of coordinates $x+iy \rightarrow x$ and $x-iy \rightarrow y$ the equation of the surface in the new coordinate is of the following form: 

\begin{align*}
z^2 = P(x) Q(y),
\end{align*}
where $P(x)= (x- z_1) \cdots (x-z_m)$ , $Q(x) = (x-\overline{z}_1) \cdots (x- \overline{z}_m)$ and $z_j= \alpha_j + i \beta_j$ for $j= 1 , \cdots , m $.
Now for singularities we have to compute the partial differentials yielding the following equations:

\begin{align*}
z=0 , ~~P'(x)Q(y)=0 ~~, P(x)Q'(y)=0 .
\end{align*}
Since these points should be on the surface and the fact that the  polynomials $P, Q$ has no multiple roots we should have 

\begin{align*}
z=0 , P(x)=0 , Q(y)=0 .
\end{align*}
So the singular points on the affine part ($w=1$) are the $m^2$ points  $(z_i ,\overline{ z}_j , 0)$ for $i,j =1 , \ldots m$.  One can easily see that the point at infinity is a singular point of the following type:  
$$z^{2m-2}= x^m y^m , $$
Also the type of singularity at singular points $(z_i ,\overline{ z}_j , 0)$ is of the type $z^2 = x y$.
In the next step we show that the singular surface $S$ with equation $ z^2 = P(x) Q(y)$ ($P$ and $Q$ are two polynomials of the same degree without multiple roots) is a surface of general type.

\medskip

In fact the surface $X$ is the branched double cover of the surface $\mathbb{P}^1 \times \mathbb{P}^1$  by the morphism  $\pi: (x, y, z) \rightarrow (\frac{x}{z}, \frac{y}{z})$ which is clearly a rational morphism.  This morphism is clearly branched along the locus $z=0$ or equivalently  $P(x) Q(y) =0$, that is the union of the $2g+2$ fibers from each ruling of  $\mathbb{P}^1 \times \mathbb{P}^1$.\\
Now we can see that the surface $X$ is a surface of general type.  Hence we showed that the surface $X$ is a singular surface with singularity at the $m^2 =(2g+2)^2$ points $(z_i ,\overline{ z}_j , 0)$
lying over the double points of the branched divisor.  
We can simply write down the canonical differentials on $X$: if we set 
\begin{align*}
\omega : = \frac{dy \wedge dx}{z},
\end{align*}
we can readily check that the products
\begin{align*}
\omega_{k, l}:  =  y^k x^l \omega
\end{align*}

give regular canonical differentials on the smooth locus of the surface $X$ for all $0\le k, l \le g-1$.   One can easily see that the ramification divisor $R$ of $\pi$ is the zero divisor of the function $z$ and we have the equality

\begin{align*}
2R =(z^2)= (P(x)) + (Q(y)),
\end{align*}

Applying the Riemann-Hurwitz formula to the cover $\pi : S \rightarrow \mathbb{P}^1 \times \mathbb{P}^1$ on the smooth part of $X$ we obtain:


\begin{align*}
K_X &= \pi^* (K_{\mathbb{P}^1 \times \mathbb{P}^1 }) \otimes \mathcal{O}_X(R) \\
&=   \pi^* (K_{\mathbb{P}^1 \times \mathbb{P}^1 }) \otimes \mathcal{O}_{\mathbb{P}^1 \times \mathbb{P}^1}(g+1 , g+1) \\
&= \pi^* (  \mathcal{O}_{\mathbb{P}^1 \times \mathbb{P}^1}(g-1 , g-1)) , 
\end{align*}
which is ample.

Now we need to deal with singular points of the surface $X$.  At the singular point $( z_i ,\overline{ z}_j , 0)$ after replacing $x$ and $y$ by $x- z_i$ and $y - \overline{ z}_j  $ respectively, the surface $X$ will have the local equation $z^2 = x y$.
Hence by blowing up once we obtain a smooth surface $\tilde{X}$.  Moreover, in terms of local coordinates we have
\begin{align*}
z' = z , ~ ~ x' = x/z , ~~ y' = y/z
\end{align*}
on $\tilde{X}$.  The pullback of the form $\omega$ to the surface $\tilde{X}$ may be written as 

\begin{align*} 
\omega' &= \frac{d(y' z') \wedge d(x' z')}{z'}\\
&=z' ( dy' \wedge dx' ) +x' (dy' \wedge dz') + y' (dz' \wedge dx') ,
\end{align*}
which is regular on all of $\tilde{X}$.  Thus all the forms $\omega_{k, l}$ above pull back to regular forms on $\tilde{X}$.  Therefore in this way we proved that the Kodaira dimension of the surface $\tilde{X}$ (and also for the singular surface $X$)  is equal to $2$.   This complete the proof of Theorem 1. $\blacksquare$

\medskip

Following the above theorem, we can pose the similar plausible question in $\mathbb{R}^n$ rather than $\mathbb{R}^2$.  But the method of the proof for the case $\mathbb{R}^2$ can't be applied to this more generalized case.

\medskip

{\bf Question}-Let $A= \{A_1 , \ldots , A_m\} $ be a finite subset of points in the $n$-dimensional Euclidean space not all on a hyperplane.  Then can one deduce that the associated distance hypersurface (its projectivization)
\begin{equation}
x_{n+1}^2=\prod_{i=1}^m ( (x_1 - a_{i1})^2 + (x_2 - a_{i2})^2 + \dots + (x_n - a_{in})^2 )
\end{equation}

 in $\mathbb{P}^{n+1}$ is of general type?\\

\begin{rem}
By the above theorem we can in fact  generate surfaces of general type of arbitrary degree $d \ge 12$ in the projective space $\mathbb{P}^3$.  

\end{rem}

Now we can prove Theorem \ref{ulam} which answers the question asked by Ulam and Erd\'os.\\

{\bf Proof of Theorem 2}~-
We  apply Theorem \ref{gentype} to resolve Ulam conjecture.  Assume $S \subset \mathbb{R}^2$ is a dense rational distance subset of the plane and also assume that $(0,0), (1,0)\in S$.  By lemma \ref{a0} there exists a square free integer $k$ for which the points in $S$ are of the form $(a,b\sqrt{k})$, where $a, b \in \mathbb{Q}$.  Since $S$ is dense we can choose  $6$ points $A=\{(a_1,b_1\sqrt{k}), \cdots, (a_6, b_6 \sqrt{k})\}$ in general position in $S$.   Let $X$ be the distance surface associated to the subset $A$ which is given by the equation 
\begin{equation}
z^2 =  \prod_{i=1}^{6}\{(x- a_i)^2 + (y- b_i\sqrt{k})^2\},
\end{equation}
and it is clearly defined over the number field $K= \mathbb{Q}(\sqrt{k})$.  By the Theorem \ref{gentype} we know that the surface $X$ defined over $K$ is a surface of general type and assuming  Bombieri-Lang conjecture the $K$-rational points on $X$ can not be Zariski dense.  Now since $S$ is dense in $\mathbb{R}^2$ it is in fact Zariski dense in $\mathbb{C}^2$ and so is dense in $\mathbb{P}^1 \times \mathbb{P}^1$.  On the other hand the surface $X$ is birationally the double cover of the surface $\mathbb{P}^1 \times \mathbb{P}^1$ by the map $(x, y, z) \rightarrow (\frac{x}{z}, \frac{y}{z})$, and so the dense set $S$ ￼ generates a Zariski-dense set of $K$-rational points in the surface $X$ which contradicts the Bombieri-Lang conjecture. $\blacksquare$ 
\\

{\bf Proof of Corollary \ref{cor1}} - As we see in the proof of Theorem 2 the set $S$ is not Zariski dense in $\mathbb{R}^2$.  Hence $S$ should be contained in the union of finitely many irreducible algebraic curves $C_1, \ldots , C_n$ in the plane.  According to Theorem 3 lines and circles are the only irreducible algebraic curves that contain an infinite rational set, so one of the irreducible curves $C_i$ should be a line (circle) with infinitely many points of rational set $S$ and  according Theorem 4 all but $4$ (resp. $3$) points of $S$ are on the line (resp. circle). $\blacksquare$\\

\medskip

I believe that if the points in the subset $A=\{ (\alpha_1, \beta_1), \cdots, (\alpha_n, \beta_n)\}$ are chosen generically then the rational curves in the corresponding distance surface $S: z^2 = P(x) Q(y)$ are contained in the subvariety $D=\{z=0\} \subset S$.  This problem is interesting enough to be dealt with in a future research.

\medskip
In the final step of preparation of this paper I noticed that Terence Tao has given another conditional proof of this problem in his weblog which is more difficult than mine.  In fact his example for a surface of general type is a complete intersection of $4$ hyper-surfaces in $\mathbb{C}^6$.\\

{\bf Aknowledgement}- I would like to thank Professor Tao for giving me some  instructive comments that improved the original version of this paper and suggesting me to publish this result.  I am also grateful to M. Makhul who  brought to my attention this problem and many stimulating discussions.  I also thank K. Shokri for reading the first draft and useful suggestions. The author was supported by a grant from Shahid Beheshti University. 


\bigskip

~\\
[3mm]
\noindent{ {\scriptsize \textsc {Department of Mathematical Sciences, Shahid Beheshti University, G.C., P.O.Box  {\rm 1983963113}, Tehran, Iran.}}} \\
 [1mm]
\noindent{ {\scriptsize \textsc {Email: $J_{-}Shaffaf@sbu.ac.ir$}}}\\


\end{document}